\documentclass[a4paper]{article}

\usepackage{amssymb,amsmath,amsthm}
\usepackage[T1]{fontenc}
\usepackage[utf8]{inputenc}
\usepackage{url}
\usepackage{enumitem}
\usepackage{color}

\newcommand{\Q}{\mathbb{Q}}

\theoremstyle{plain}
\newtheorem{theorem}{Theorem} 

\newtheorem*{proposition*}{Proposition}

\newtheorem{conjecture}[theorem]{Conjecture}
\newtheorem*{conjecture*}{Conjecture}

\theoremstyle{definition}

\theoremstyle{remark}

\newtheorem*{remark*}{Remark}

\newcommand*\patchAmsMathEnvironmentForLineno[1]{%
  \expandafter\let\csname old#1\expandafter\endcsname\csname #1\endcsname
  \expandafter\let\csname oldend#1\expandafter\endcsname\csname end#1\endcsname
  \renewenvironment{#1}%
     {\linenomath\csname old#1\endcsname}%
     {\csname oldend#1\endcsname\endlinenomath}}%
\newcommand*\patchBothAmsMathEnvironmentsForLineno[1]{%
  \patchAmsMathEnvironmentForLineno{#1}%
  \patchAmsMathEnvironmentForLineno{#1*}}%
\AtBeginDocument{%
\patchBothAmsMathEnvironmentsForLineno{equation}%
\patchBothAmsMathEnvironmentsForLineno{align}%
\patchBothAmsMathEnvironmentsForLineno{flalign}%
\patchBothAmsMathEnvironmentsForLineno{alignat}%
\patchBothAmsMathEnvironmentsForLineno{gather}%
\patchBothAmsMathEnvironmentsForLineno{multline}%
}
\usepackage{lineno}

\title{Counting pop-stacked permutations in polynomial time}

\author{
		Anders Claesson\\
		\small Division of Mathematics\\
		\small The Science Institute\\
		\small University of Iceland\\
		\small \texttt{akc@hi.is}
	\and
		Bjarki \'Ag\'ust Gu\dh{}mundsson\\
		\small Department of Computer Science\\
		\small Reykjav\'ik University\\
		\small Reykjav\'ik, Iceland\\
		\small \texttt{bjarkig@ru.is}
	\and
		Jay Pantone\\
		\small Department of Mathematical and Statistical Sciences\\
		\small Marquette University\\
		\small Milwaukee, WI, USA\\
		\small \texttt{jay.pantone@marquette.edu}
}

\date{}

\begin{document}
\maketitle
\thispagestyle{empty}

\begin{abstract}
  Permutations in the image of the pop-stack operator are said to be
  \emph{pop-stacked}. We give a polynomial-time algorithm to count pop-stacked
  permutations up to a fixed length and we use it to compute the first 1000
  terms of the corresponding counting sequence. Only the first 16
  terms had previously been computed.  With the 1000 terms we prove some
  negative results concerning the nature of the generating function for
  pop-stacked permutations. We also predict the asymptotic behavior of
  the counting sequence using differential approximation.
\end{abstract}

\section{Pop-stacked permutations}

The abstract data type known as a \emph{stack} has two operations:
\emph{push} adds an element at the top of the stack; \emph{pop} removes
the top element from the stack. A \emph{pop-stack} is a variation of
this introduced by Avis and Newborn~\cite{avis1981pop} in which the pop
operation empties the entire stack.

Let $\pi = a_1a_2 \dots a_n$ be a permutation of
$[n]=\{1,2,\dots,n\}$. An \emph{ascending run} of $\pi$ is a maximal
sequence of consecutive ascending letters
$a_i < a_{i+1} < \dots < a_{i+d-1}$, and
a \emph{descending run} is defined similarly.  For instance, the
ascending runs of $\pi = 617849235$ are $6$, $178$, $49$ and $235$; its
descending runs are $61$, $7$, $84$, $92$, $3$ and $5$.

Let $P(\pi)$ be the result of greedily sorting $\pi$ using a pop-stack
subject to the constraint that elements on the pop-stack are increasing
when read from the top to the bottom of the stack. In other words, if we
factor $\pi$ into its descending runs $\pi = D_1 D_2 \dots D_m$, then
$P(\pi)$ is obtained by reversing each of those runs:
$P(\pi)=D_1^r D_2^r \dots D_m^r$. For instance, $P(5321764)=1235467$ and
$P(617849235) = 167482935$. A permutation $\pi$ is said to be sortable
by a pop-stack if $P(\pi)$ is the identity permutation. More generally,
$\pi$ is said to be sortable by $k$ passes through a pop-stack if $P^k(\pi)$
is the identity permutation. Claesson and Guðmundsson~\cite{ClGu2019}
showed that the generating function for the number of permutations of
$[n]$ that are sortable by $k$ passes through a pop-stack is always
rational.

Asinowski et~al.~\cite{asinowski-etal-2019} defined that $\sigma$ is
\emph{pop-stacked} if $\sigma=P(\pi)$ for some permutation $\pi$, and
gave the following theorem.

\begin{theorem}[Asinowski et~al.~\cite{asinowski-etal-2019}]
  A permutation is pop-stacked if and only if for each pair
  $(R_i,R_{i+1})$ of its adjacent ascending runs
  $\min R_i < \max R_{i+1}$.
\end{theorem}

They further showed that the generating function for pop-stacked
permutations of $[n]$ with exactly $k$ ascending runs is rational for
each $k$. Enumerating pop-stacked permutations without this restriction
is, however, an open problem. Asinowski et~al.\ initiated an investigation into
this by calculating the number of pop-stacked permutations of length
$n=1,\ldots,16$, adding the resulting sequence to the OEIS~\cite{oeis-no-link}
as A307030 and noting that ``this sequence is hard to compute''.
In the following section, we give an efficient algorithm for counting pop-stacked
permutations, expanding the sequence up to $n=1000$. While the algorithm and
the augmented sequence could give additional insight into the structure of
pop-stacked permutations, finding a generating function or a closed form
solution to their enumeration remains an open problem. Section~\ref{section:experimental} gives experimental data in this direction.

\section{Polynomial-time counting algorithm}

A \emph{ballot}, alternatively known as an ordered set partition, is a
collection of pairwise disjoint nonempty sets, referred to as blocks,
where the blocks are assigned some total ordering. Any permutation can
be seen as a ballot by decomposing it into its ascending runs. The permutation
$\pi = 6 178 49 235$ would then be viewed as the ballot
$\{6\}\{1,7,8\}\{4,9\}\{2,3,5\}$. Conversely, a ballot
$B_1 B_2\dots B_k$ represents a permutation in this manner if,
and only if, $\max B_i > \min B_{i+1}$ for each $i$ in $[k-1]$. Thus,
the ballots corresponding to pop-stacked permutations are precisely
those such that
\[
	\max B_i > \min B_{i+1} \quad\text{and}\quad
	\min B_i < \max B_{i+1}.
\]
In other words, the intervals between the smallest and largest elements
of each pair of adjacent blocks overlap,
$$[\min B_i, \max B_i] \cap [\min B_{i+1}, \max B_{i+1}]\neq\emptyset,
$$
and we call these ballots \emph{overlapping}; here, $[a,b]$ denotes the
interval $\{a,a+1,\dots,b\}$.
Let $F[U]$ be the set
of overlapping ballots whose underlying set is $U$. As an example,
\[
	F[\{1,2,3\}] = \bigl\{\{1, 2, 3\}, \{2\}\{1, 3\}, \{1, 3\}\{2\}\bigr\}.
\]
Let $F_{c,d}[U]$ denote the subset of $F[U]$
whose last block, $B$, is such that $c = \min B$ and $d = \max B$.
Clearly, if $c > d$ then $F_{c,d}[U]=\emptyset$. Also,
\[
	F[U] \,=\!\bigcup_{c,d\,\in\,U}\! F_{c,d}[U].
\]

If $c = \min U$ and $d = \max U$, then one possibility is that there is
a single block consisting of all elements of $U$. Let us now consider
the more typical case when there are two or more blocks, and let us write
the ballot as $B_1 B_2\dots B_k$. By definition, its last block, $B_k$,
satisfies $c = \min B_k$ and $d = \max B_k$, or expressed differently
$\{c,d\}\subseteq B_k \subseteq [c,d]$. Let $a = \min B_{k-1}$ and
$b = \max B_{k-1}$. The blocks $B_{k-1}$ and $B_k$ overlap if, and only
if, $a < d$ and $b > c$. Thus
\begin{equation}\label{rec-decomp}
\begin{aligned}
F_{c,d}[U]
\;=\;\;& \{ U : c = \min U \land d = \max U \}
\;\;\cup \\
&\bigcup_{
  \substack{
    \{c,d\}\,\subseteq\, B\,\subseteq\, [c,d]\smallskip\\
    a,b \,\in\, U\setminus B \smallskip \\
    a < d \,\land\, b > c
  }
}
F_{a,b}[U \setminus B]\, B,\qquad\qquad
\end{aligned}
\end{equation}
where $F_{a,b}[U \setminus B]B$ is the set
$\bigl\{ wB : w \in F_{a,b}[U\setminus B]\bigr\}$, and the somewhat
cryptic looking $\{ U : c = \min U \land d = \max U \}$ expresses the
singleton $\{U\}$ if $c = \min U$ and $d = \max U$, and the empty set
otherwise.

We now turn to counting. Let $f(n)$ be the number of overlapping ballots
of $[n]$. That is, $f(n) = |F[n]|$ in which $F[n]$ is
short for $F[\{1,\dots,n\}]$.
Also,
let $f_{c,d}(n) = |F_{c,d}[n]|$. If $c > d$ then $f_{c,d}(n) =
0$. Otherwise we shall use the recursive decomposition~\eqref{rec-decomp} and
do case analysis based on whether $c$ and $d$ are the same or two
distinct elements.

If $c=d$, then the last block consists of a single point. In terms of
\eqref{rec-decomp} the ballot is written $w\{c\}$, where
$w\in F_{a,b}[[n]\setminus\{c\}]$ and $a<c<b$. After ``rescaling'' we
can consider $w$ a ballot in $F[n-1]$; here we subtract $1$ from each
element larger than $c$. Note that this, however, also lowers the value
of $b$ by one. Thus, the number of such ballots is
\[
	\sum_{a=1}^{c-1} \sum_{b=c}^{n-1} f_{a,b}(n-1).
\]

If $c < d$, then write the ballot as $wB$ and let $\ell=|B|-2$. There
are $\binom{d-c-1}{\ell}$ ways to choose $B$. After rescaling we have
$w\in F_{a,b}[n-\ell-2]$, where $a \leq d-\ell-2$ and $b\geq c$.
Thus, the number of such ballots is
\[
	\sum_{\ell=0}^{d-c-1} \binom{d-c-1}{\ell} \sum_{a=1}^{d-\ell-2} \sum_{b=c}^{n-\ell-2} f_{a,b}(n-\ell-2).
\]

Finally, if $c=1$ and $d=n$, we also count the case where the ballot consists
of a single block. Taking all this together, we have that
\begin{equation} \label{eq:f-recurrence}
\begin{split}
    f_{c,d}(n)
    &= [c = 1 \land d=n] \\
    &+ [c=d]\, \sum_{a=1}^{c-1} \sum_{b=c}^n f_{a,b}(n-1) \\
    &+ [c<d]\, \sum_{\ell=0}^{d-c-1} \binom{d-c-1}{\ell} \sum_{a=1}^{d-\ell-2} \sum_{b=c}^{n-\ell-2} f_{a,b}(n-\ell-2).
\end{split}
\end{equation}
Here $[p]$ is the Iverson bracket: it converts the proposition $p$ into
$1$ if $p$ is satisfied, and $0$ otherwise. Further,
$f(n) = \sum_{a=1}^n \sum_{b=a}^n f_{a,b}(n)$.

Recurrence~\eqref{eq:f-recurrence} can be augmented to count overlapping
ballots with a specific number of blocks, or, equivalently, pop-stacked
permutations with a specific number of ascending runs. Let $f_{c,d}(n,k)$ denote
the number of overlapping ballots of $[n]$ with exactly $k$
blocks. Then we have $f(n,k) = \sum_{a=1}^n \sum_{b=a}^n f_{a,b}(n,k)$ and
\begin{equation} \label{eq:f-k-recurrence}
\begin{split}
  f_{c,d}(n,k)
  &= [c = 1 \land d=n \land k=1] \\
  &+ [c=d]\, \sum_{a=1}^{c-1} \sum_{b=c}^n f_{a,b}(n-1,k-1) \\
  &+ [c<d]\, \sum_{\ell=0}^{d-c-1}\binom{d-c-1}{\ell}\sum_{a=1}^{d-\ell-2}\sum_{b=c}^{n-\ell-2}f_{a,b}(n-\ell-2,k-1).
\end{split}
\end{equation}

Note that there are two locations in the recurrence
\eqref{eq:f-recurrence} where we have a plain two-dimensional sum over
$f$, that is
$\sum_{a=\star}^{\star}\sum_{b=\star}^{\star} f_{a,b}(\star)$, where
$\star$ are fixed and not dependent on $a$, $b$ or each other. We
simplify these two-dimensional sums using ``prefix sums''. Let
\[
    g_{c,d}(n) = \sum_{a=1}^c \sum_{b=1}^d f_{a,b}(n)
\]

In particular, $g_{c,d}(n) = 0$ if $c = 0$ or $d=0$. Note that
\begin{equation}\label{eq:g_recurrence_fast}
    g_{c,d}(n) = f_{c,d}(n) + g_{c-1,d}(n) + g_{c,d-1}(n) - g_{c-1,d-1}(n).
\end{equation}

Also noting that
\[
	\sum_{a=p}^q \sum_{b=r}^s f_{a,b}(n) = g_{q,s}(n) - g_{p-1,s}(n) - g_{q,r-1}(n) + g_{p-1,r-1}(n),
\]
we can now simplify the above equation to
\begin{equation}\label{eq:f_recurrence_fast}
\begin{split}
  f_{c,d}(n)
  &= [c = 1 \land d = n] \\[1.2ex]
  &+ [c=d]\, \Delta_{c-1,n,c-1}(n-1) \\[1ex]
  &+ [c<d] \sum_{\ell=0}^{d-c-1} \binom{d-c-1}{\ell}\Delta_{d-2-\ell,n-2-\ell,c-1}(n-2-\ell)
\end{split}
\end{equation}
where $\Delta_{u,v,w}(n) = g_{u,v}(n) - g_{u,w}(n)$. We further have
$f(n) = g_{n,n}(n)$. The same simplification can also be applied to the
recurrence for counting by blocks.


Say we wanted to compute $f(n)$ for all $1\leq n\leq N$. We can precompute
binomial coefficients $\binom{n}{k}$ for all $0 \leq k \leq n \leq N$ using the
recurrence $\binom{n}{k} = \binom{n-1}{k-1} + \binom{n-1}{k}$. Then, using
dynamic programming we can compute $f_{c,d}(n)$, $g_{c,d}(n)$ and $f(n)$ using
Recurrences~\ref{eq:g_recurrence_fast} and \ref{eq:f_recurrence_fast} for all
$1 \leq c$, $d \leq n \leq N$ in $O(N^4)$ time using $O(N^3)$ memory. When
counting by blocks this is $O(N^5)$ time, but $O(N^3)$ memory is still
sufficient.

This assumes that all arithmetic operations are $O(1)$. In reality, some of the
numbers are on the order of $N!$. This means that multiprecision arithmetic has
to be used, which slows down the computation considerably. One way to speed this up
is to choose a set of relatively small primes whose product is greater than $N!$.
For each prime $p$, the above computation is then carried out in the finite
field $\textbf{F}_p$. This can be done in parallel, as the computation for
different primes is independent. The values of $f(n)$, which are guaranteed to be at
most $N!$ for all $n\leq N$, are then recovered using the Chinese Remainder Theorem.

This was used to calculate the number of pop-stacked permutations of each
length up to $N=1000$. With $286$ distinct primes just under $10^9$, and one
CPU core per prime, the computation took just under an hour to complete, with
each core using 3.8GiB of RAM. In a similar manner the number of pop-stacked
permutations of each length up to $N=300$ grouped by number of ascending runs
were computed. Table~\ref{tbl:seq} gives the number of pop-stacked permutations
of each length up to $N=45$, but the complete results, along with the code used
to generate the results, can be found on GitHub~\cite{github}.

\begin{table}
    \centering
    \begin{tabular}{rr}
        $n$ & $f(n)$ \\
        \hline
        1 & 1 \\
        2 & 1 \\
        3 & 3 \\
        4 & 11 \\
        5 & 49 \\
        6 & 263 \\
        7 & 1653 \\
        8 & 11877 \\
        9 & 95991 \\
        10 & 862047 \\
        11 & 8516221 \\
        12 & 91782159 \\
        13 & 1071601285 \\
        14 & 13473914281 \\
        15 & 181517350571 \\
        16 & 2608383775171 \\
        17 & 39824825088809 \\
        18 & 643813226048935 \\
        19 & 10986188094959045 \\
        20 & 197337931571468445 \\
        21 & 3721889002400665951 \\
        22 & 73539326922210382215 \\
        23 & 1519081379788242418149 \\
        24 & 32743555520207058219615 \\
        25 & 735189675389014372317381 \\
        26 & 17167470189102029106503457 \\
        27 & 416297325393961581614919699 \\
        28 & 10468759109047048511785181499 \\
        29 & 272663345523662949571086535201 \\
        30 & 7346518362495550669587951987399 \\
        31 & 204539324291355079758576427320853 \\
        32 & 5878416448467628215599958670190869 \\
        33 & 174223945386975482728912851110751431 \\
        34 & 5320106374135453888563313157982976111 \\
        35 & 167232974698164950641578719412434688845 \\
        36 & 5407019929661274797886581276653666104943 \\
        37 & 179677314965899717327756420597568210468933 \\
        38 & 6132116544121046402686046213590718114272089 \\
        39 & 214787281796488809444762543177377466419782267 \\
        40 & 7716175695131570964771559074490172330993576115 \\
        41 & 284131588386675257705011846785657928372695002841 \\
        42 & 10717718945463416620327720805595647805635809236711 \\
        43 & 413908527884993695909526722330319436067536797304549 \\
        44 & 16356508568742954048255540186930772843919017766669517 \\
        45 & 661053598808034620660440013405109251647269697650963759
    \end{tabular}
    \caption{The number of pop-stacked permutations of each length up to $N=45$.}
    \label{tbl:seq}
\end{table}

\section{Experimental analysis}
\label{section:experimental}

With the first 1000 terms of the counting sequence of pop-stacked permutations now calculated, we turn to a pair of experimental techniques for an empirical analysis: \emph{automated fitting} and \emph{differential approximation}. Given initial terms of a counting sequence, the first of these methods searches for a generating function whose power series expansion matches the sequence, while the second predicts the asymptotic growth of the sequence.

For the counting sequence at hand, automated fitting does not conjecture a generating function, giving instead several (rigorous) negative results, while differential approximation gives very precise estimates of the asymptotic behavior.

\subsection{Automated fitting for pop-stacked permutations}

Let $a_0,a_1,\ldots$ be a counting sequence and $F(x) = \sum_{n \geq 0}a_nx^n$ its generating function. If $F(x)$ is a rational function, then we can write $F(x) = p(x)/q(x)$ for relatively prime polynomials $p(x), q(x) \in \Q[x]$; equivalently,
\begin{equation}
	\label{equation:rational-form}
	q(x)F(x) - p(x) = 0.
\end{equation}
\sloppy Conversely, suppose we are given only some initial terms $a_0, a_1, \ldots, a_n$ of a counting sequence and want to determine whether the generating function $F(x)$ of the unknown counting sequence is rational. If $F(x)$ is rational with $\max(\deg(p(x)),\deg(q(x))) = d$, then we can write Equation~(\ref{equation:rational-form}) as
\begin{equation}
	\label{equation:generic-rational}
	(q_0 + q_1x + \cdots + q_dx^d)(a_0 + a_1x + \cdots + a_nx^n) - (p_0 + p_1x + \cdots + p_dx^d) = 0.
\end{equation}
Expanding the left-hand side gives a polynomial in $x$, and the coefficients of $x^0, x^1, \ldots, x^n$ must all equal $0$. We thus have a system of $n+1$ equations in the $2d+2$ unknowns $p_0,\ldots,p_d,q_0,\ldots,q_d$. A generic system of this form is likely to have non-trivial solutions when $n \leq 2d$, and so when initial terms up to $a_n$ are known, it is only productive to consider $d$ such that $2d < n$.

If this system has no non-trivial solution, then we are guaranteed that $F(x)$ is not rational  with numerator and denominator of degree at most $d$. If the system does have a non-trivial solution, then it is possible, though far from guaranteed, that
\[
	F(x) = \frac{p_0 + p_1x + \cdots + p_dx^d}{q_0 + q_1x + \cdots + q_dx^d}.
\]
The larger the difference between $n$ and $2d$, the more confident that one can be in such a conjecture. Empirically, this is like using the first $2d$ known terms to guess the rational generating function and the remaining $n-2d$ as confirmation.

Automated fitting can be extended beyond the realm of rational generating functions. A generating function $F(x)$ is called \emph{algebraic} if there are polynomials $p_0(x), \ldots, p_m(x) \in \Q[x]$ such that
\[
	p_m(x)F^m(x) + \cdots + p_1(x)F(x) + p_0(x) = 0,
\]
called \emph{differentially finite} (or \emph{D-finite}) if there are polynomials $p_0(x), \ldots, p_k(x), q(x) \in \Q[x]$ such that
\[
	p_k(x)F^{(k)}(x) + \cdots + p_1(x)F'(x) + p_0(x)F(x) + q(x) = 0,
\]
and called \emph{differentially algebraic} (or \emph{D-algebraic}) if there exists a $(k+2)$-variate polynomial $P$ with coefficients in $\Q$ such that
\[
	P(x, F(x), F'(x), \ldots, F^{(k)}(x)) = 0.
\]

To determine whether a generating function $F(x)$ is algebraic given some initial terms, an equation similar to~(\ref{equation:generic-rational}) can be set up assuming each $p_i(x)$ has degree at most $d$, giving a linear system with $n$ equations and $(m+1)(d+1)$ unknowns. In the D-finite case, the system has $(k+2)(d+1)$ unknowns. The D-algebraic case requires further assumptions about form---the ideas are similar, but not worth elaborating upon here. There are various software packages that perform fitting of this kind, including \texttt{Gfun}~\cite{salvy:gfun} in Maple, \texttt{Guess}~\cite{kauers:guess} in Mathematica, and \texttt{Guess}~\cite{rubey:fricas-guess} in FriCAS. We have used a different package, \texttt{GuessFunc}, written by the third author.

We applied automated fitting to the counting sequence of pop-stacked permutations up to length $1000$, and found no conjectured rational, algebraic, D-finite, or D-algebraic form for the unknown generating function $F(x)$. From this we can conclude rigorously that, for example,
\begin{enumerate}[label={$\diamond$}]
	\item If $F(x)$ is rational, then either the degree of the denominator or the degree of the numerator is at least $500$.
	\item If $F(x)$ is algebraic, then the degree of algebraicity $m$ and the maximum degree of polynomial coefficient $d = \max(p_0(x), \ldots, p_m(x))$ must satisfy $(m+1)(d+1) > 1000$.
	\item If $F(x)$ is D-finite, then the differential order $k$ and the maximum degree of polynomial coefficient $d = \max(q(x), p_0(x), \ldots, p_k(x))$ must satisfy $(k+2)(d+1) > 1000$.
\end{enumerate}
A similar negative result could be written for the D-algebraic case, although it would require further explanation of the structure of the corresponding search space.

One can also apply various transformations to the generating function before initiating the automated fitting procedure. In addition to trying to find a fit for the ordinary generating function $F(x) = \sum_{n \geq 0}a_nx^n$, we also attempted to find a fit for the exponential generating function $\sum_{n \geq 0} (a_n/n!)x^n$, the reciprocal $1/F(x)$, the compositional inverse $F(x)^{\langle -1 \rangle}$, and also several combinations of these transformations. No results were found.

\subsection{Automated fitting for pop-stacked permutations with a fixed number of ascending runs}

Let $F_k(x)$ denote the power series for those pop-stacked permutations with
precisely $k$ ascending runs. Asinowski et~al.~\cite{asinowski-etal-2019}
showed that these permutations are in bijection with words from a
regular language that is recognized by a certain deterministic finite automaton
(DFA) $\mathcal{A}_k$, proving that $F_k(x)$ is rational. Furthermore, a system
of linear equations can be derived from this DFA, whose solution gives
$F_k(x)$. Deriving $F_k(x)$ in this way is only practical for small values of
$k$, however, as the number of states in $\mathcal{A}_k$ grows exponentially
with $k$.

As mentioned earlier, Recurrence~(\ref{eq:f-k-recurrence}) permits the fast
computation of the counting sequence for pop-stacked permutations with a fixed
number of ascending runs. This, along with the techniques of automated fitting
gives rise to a different approach for finding $F_k(x)$, albeit
heuristically\footnote{Given enough terms of the sequence, automated fitting
will find $F_k(x)$. The number of terms required is the sum of the degrees of the numerator and denominator of $F_k(x)$,
which is not known. An upper bound is twice the number of states in
$\mathcal{A}_k$, which is exponential.}.

Using the counting sequence for pop-stacked permutations of length at most
$300$ with a fixed number of ascending runs, we were able to find a rational
fit for each $F_k(x)$ for $k \leq 24$. We were further able to verify that the
rational fits were exact for $k\leq 6$ by using the previously mentioned method
based on Asinowski et~al.~\cite{asinowski-etal-2019}. The first four generating
functions follow.
\begin{align*}
	F_1(x) &= \frac{x}{1-x},\\
	F_2(x) &= \frac{2x^3}{(1-2x)(1-x)^2},\\
	F_3(x) &= \frac{2x^4(1+3x-6x^2)}{(1-3x)(1-2x)^2(1-3x)^3},\\
	F_4(x) &= \frac{2x^6(21-74x+5x^2+180x^3-144x^4)}{(1-4x)(1-3x)^2(1-2x)^3(1-x)^4}.
\end{align*}

Based on this data, which can be found in full on GitHub~\cite{github}, we pose
the following conjecture.
\begin{conjecture}
	For all $k$, the rational generating function $F_k(x)$ can be written as
	\[
		F_k(x) = N_k(x)\Big/\prod_{i=1}^k (1-ix)^{k-i+1},
	\]
	where $N_k(x)$ is a polynomial of degree $k(k + 1)/2$, the same degree as the conjectured denominator.
\end{conjecture}

\subsection{Differential approximation}

Differential approximation empirically estimates the asymptotic growth of a counting sequence based on its initial terms by using linear differential equations to model the unknown generating function and studying the complex singularities of solutions of those linear differential equations. Here we will only present the results of this analysis---for information about how differential approximation works we refer the reader to~\cite{guttmann:asymptotic-analysis, guttmann:series-analysis}.

The cornerstone of analytic combinatorics is the observation that the asymptotic behavior of a counting sequence is intimately connected to the singularities of its generating function when treated as a complex function. For example, the location of the singularities closest to the origin (the \emph{dominant} singularities) roughly determine the exponential growth of the counting sequence, and the nature of those singularities determines the sub-exponential behavior. 

The output of differential approximation is an estimate of the location and nature (specifically, the critical exponent) of all singularities of the unknown generating function based on the given known initial terms. Typically, although not always, the dominant singularity is predicted with the highest precision, with the precision of the estimates of other singularities decreasing as distance from the origin increases. Obviously such an analysis is only experimental, but in practice the estimates given by differential approximation are incredibly accurate. In tests where the true singularity structure of a generating function is independently known, the estimates from differential approximation are rarely off by more than the last decimal place.

The counting sequence of pop-stacked permutations grows superexponentially~\cite{asinowski-etal-2019}, implying that its generating function has a singularity at the origin. Accordingly, we use differential approximation to analyze the exponential generating function. It predicts a number of singularities on the positive real axis, located at the values below.

\begin{footnotesize}
\begin{align*}
		&1.113439041736727043761661526918083240141390165833449466152700785053219911270\ldots\\
		&2.417184228722564007388473547672885752580057534770845001690528350200102151036\ldots\\
		&3.076673197412146436807595671137309181422151285506943038305240180949212077913\ldots\\
		&3.527590791728018755531106354662725269743465863978439496914729951030934478987\ldots\\
		&3.872438162423457670453537298789680569472671309363632792004917259462379566078\ldots\\
		&4.152519207830100565666605055176411745894938982832118599384868016797119166567\ldots\\
		&4.388766437824164163366758081274636520883940965171626205159043874261749420137\ldots\\
		&4.593300493040369902037314403433340137408669134838327397901215132095535249496\ldots\\
		&4.773787732301263733990448984231076188826829730174328444872240429327757789160\ldots\\
		&4.93539355029443080528699130532727322201728351298582403913\\
		&5.08176797057144544489527338196678922218609719159\\
		&5.215588012778242472294262722856995906\\
		&5.453200964209036692\\
		&5.55979961612\\
		&5.659669
\end{align*}
\end{footnotesize}

Each of these singularities is predicted to have critical exponent $-1$, making them simple poles. The topmost 9 estimates have been truncated to fit on the page. In reality, they are given to many more decimal places---nearly $800$ for the dominant singularity. More precise estimates could be obtained if desired. These results suggest that the exponential generating function may posses an infinite number of singularities. If true, this would imply the non-D-finiteness of both the ordinary and exponential generating functions.

Differential approximation also predicts several complex pairs of singularities, also simple poles, of which we'll list just a few.

\begin{footnotesize}
\begin{align*}
	&0.4279380975440727242991591373540946029637854497521857134254777354059489934\ldots\\
	& \, \pm 3.6012595134274782137294551323567899146878282109407492350988015900552787045\ldots i\\
	&1.8079319224525533045652715650438553186508451786578693412247786970810774117\ldots\\
	& \, \pm 4.0462349876106887702897457441128645763490304850344195743880592871046130995\ldots i\\
	&2.5083998717369662727687249193314945476381464747880461769920884622874845896\ldots\\
	& \, \pm 4.2416800160392329291940969204250545140382149982272394213372595306429864967\ldots i
\end{align*}
\end{footnotesize}

The dominant pole at $\mu \approx 1.11343904$ implies that the exponential growth rate of the counting sequence is 
\[
	\mu^{-1} \approx 0.8981183185746869695116759646856448\ldots,
\]
implying that the asymptotic behavior of the number of pop-stacked permutations is
\[
	a_n \sim C \cdot n! \cdot (0.898118\ldots)^n.
\]
Differential approximation does provide an estimate for the constant $C$ but this can be obtained numerically given the extremely accurate estimate for $\mu$. We find that
\[
	C \approx 0.6956885490706357679957031687241101565741983507216179232324\ldots
\]
giving the final asymptotic approximation
\[
	a_n \sim (0.695688\ldots) \cdot n! \cdot (0.898118\ldots)^n.
\]
Full decimal values for the approximated singularities and constants can also be found on GitHub~\cite{github}.

\textbf{Acknowledgements.} Computations were performed on the Garpur cluster~\cite{garpur}, a joint project between the University of Iceland and the University of Reykjavik funded by the Icelandic Centre for Research. We thank them for the use of their resources.

\bibliography{pop-stacked}{}
\bibliographystyle{plain}

\end{document}